\documentclass[11pt,a4paper]{article}
\usepackage{amsthm}
\usepackage[T1]{fontenc}
\usepackage[utf8]{inputenc}
\usepackage{lmodern}
\usepackage{amsmath,amssymb,amsthm}
\usepackage[margin=1in]{geometry}
\usepackage[hidelinks]{hyperref}
\usepackage{microtype}

\numberwithin{equation}{section}

\newtheorem{theorem}{Theorem}[section]
\newtheorem{proposition}[theorem]{Proposition}
\newtheorem{lemma}[theorem]{Lemma}
\newtheorem{corollary}[theorem]{Corollary}
\newtheorem{definition}[theorem]{Definition}
\theoremstyle{definition}
\newtheorem{example}[theorem]{Example}

\theoremstyle{remark}
\newtheorem{remark}[theorem]{Remark}

\newcommand{\R}{\mathbb{R}}
\newcommand{\Z}{\mathbb{Z}}
\newcommand{\Nscg}{N_{\mathrm{scg}}}
\newcommand{\Khat}{\widehat K}

\title{Simple closed geodesics on the Klein bottle:\\
the sharp lower bound}

\author{Gabriella Freitas \and Marcelo Miranda}
\date{}

\begin{document}

\maketitle

\begin{abstract}
We revisit the Klein bottle case of Ballmann's theorem on simple closed geodesics. The final counting argument in Ballmann's proof relies on the existence of a nontrivial homological family in a primitive one-sided free homotopy class. We identify an obstruction to this assertion and establish instead the lower bound of four geometrically distinct simple closed geodesics: two in the horizontal two-sided class and one in each of the two primitive one-sided classes. We then show that this bound is sharp. The exclusion of any additional simple closed geodesics is obtained by combining the topology of the orientation double cover with a separation of the geodesic equations and a period estimate based on Wirtinger’s inequality.
\end{abstract}

\section{Introduction}

The existence and multiplicity of closed geodesics are classical and relevant problems in Riemannian geometry. One of the first results on the existence of closed geodesics is due to Hadamard \cite{hadamard}. In 1898, he showed that a closed surface of negative curvature admits a closed geodesic. In 1917, Birkhoff established the existence of a closed geodesic on any two-dimensional sphere by means of variational methods \cite{birkhoff}. Another fundamental result was obtained by Lusternik and Schinirelmann \cite{LS1929}, who proved that every Riemannian metric on the two-dimensional sphere admits at least three simple closed geodesics.

While the existence of closed geodesics has been extensively investigated, requiring these geodesics to be simple introduces additional geometric and topological restrictions. In particular, multiplicity results for closed geodesics do not necessarily imply analogous results for simple closed geodesics. This distinction is especially significant in the study of surfaces, where the topology of embedded closed curves plays a fundamental role. Focusing on the distinction between closed geodesics and simple closed geodesics, the results of Franks and Bangert show that every Riemannian two-sphere admits infinitely many geometrically distinct closed geodesics, as can be seen in \cite{Franks1992} and \cite{Bangert1993}.

In their work, Bangert and Hingston \cite{BangertHingston1984} establish the existence of infinitely many geometrically distinct closed geodesics on compact manifolds with infinite abelian fundamental group. Furthermore, Taimanov, in \cite{Taimanov1992}, surveys interesting results on closed geodesics on surfaces. In particular, Taimanov records that every Riemannian metric on the real projective plane admits at least three geometrically distinct, noncontractible simple closed geodesics.

\begin{definition}(The Klein bottle)
 Let $K$ be the quotient of $\mathbb{R}^2$ by the group $\Gamma$ generated by the isometries
    \begin{equation}\label{eq:generators}
        \qquad
        A(x,y)=(x+2\pi,y),
        \qquad
        B(x,y)=(-x,y+2\pi),
    \end{equation}
    equipped with the induced flat metric, and let
    $\pi_K:\mathbb{R}^2\to K$ be the quotient map.
\end{definition}

Throughout the paper we use the notation
\begin{equation}\label{eq:ab-notation}
 a=[A],
 \qquad
 b=[B],
 \qquad
 \pi_1(K)=\langle a,b\mid b^{-1}ab=a^{-1}\rangle.
\end{equation}
Thus $a$ is the horizontal two-sided class, while $b$ and $ab$ are
the two primitive one-sided classes, represented by $B$ and $AB$.
The separating two-sided class is $b^2$, represented by $B^2$.

In the quotation from Ballmann in Section \ref{sec: correctquotation}, the letters retain their
original meaning: he calls the left edge $a$ and the bottom edge $b$.
For the edge identifications \eqref{eq:generators}, these correspond
to our $b$ and $a$, respectively. His product $ab$ then corresponds
to our $ba$, which is conjugate to $ab$ and hence represents the same
free homotopy class. All notation outside that quotation follows
\eqref{eq:ab-notation}.

In \cite[Theorem~2.3]{Ballmann1978}, Ballmann states that every Riemannian metric on the Klein bottle admits at least five geometrically distinct simple closed geodesics. The last paragraph of the proof obtains this number by applying a torus-type homological construction to two free homotopy classes. At least one of these classes is one-sided (see Definition \ref{def: one-two side curves}) and the asserted nontrivial $1$-cycle does not exist in the relevant quotient loop space. The following example shows that, for the flat metric, Ballmann's argument does not work:

\begin{example}\label{example: flatmetric}
    The proof of Theorem~2.3 in \cite{Ballmann1978} appears to imply
    that, for every Riemannian metric on the Klein bottle, the free
    homotopy class $ab$ contains at least two geometrically distinct
    simple closed geodesics. However, this conclusion does not hold
    for the standard flat metric, as we now show.

    We claim that
    the free homotopy class $ab$ contains exactly one simple closed
    geodesic, up to reparametrization, namely $\gamma_0(t)=\pi_K(\pi,t)$ with  $t\in\mathbb{R}/2\pi\mathbb{Z}.$

    Let $\gamma:\mathbb{R}\to K$ be a simple closed geodesic in this
    class, with prime period $T$. Since $\gamma$ represents the free
    homotopy class $ab$, we may choose a lift
    $\widetilde{\gamma}:\mathbb{R}\to\mathbb{R}^2$ such that
    \[
        \widetilde{\gamma}(T)
        =G\bigl(\widetilde{\gamma}(0)\bigr),
        \qquad
        G=A\circ B,
    \]
    where $G(x,y)=(2\pi-x,y+2\pi)$. The paths $t\mapsto\widetilde{\gamma}(t+T)$ and
    $t\mapsto G(\widetilde{\gamma}(t))$ are lifts of the same path
    $\gamma$ and agree at $t=0$. By uniqueness of path lifting,
    \[
        \widetilde{\gamma}(t+T)
        =G\bigl(\widetilde{\gamma}(t)\bigr)
        \qquad\text{for all }t\in\mathbb{R}.
    \]

    Write $\widetilde{\gamma}(t)=(x(t),y(t))$. Since $\pi_K$ is
    a local isometry, $\widetilde{\gamma}$ is a Euclidean geodesic,
    so $x'$ and $y'$ are constant. Differentiating $x(t+T)=2\pi-x(t)$ gives $x'(t+T)=-x'(t)$. As $x'$ is constant, this forces $x'\equiv0$. Substituting
    back into the preceding identity yields $x(t)\equiv\pi$.
    Hence $\gamma$ has the same image as $\gamma_0$.
\end{example}

The purpose of this paper is to correct this part of the argument and
determine the optimal lower bound. We first exhibit the obstruction
in the one-sided class and retain the valid homological construction
in the horizontal class. The resulting count is four. We then give
an explicit family of metrics for which these four geodesics are the
only simple closed geodesics, proving sharpness.

A simple closed geodesic means an embedded circle whose arclength
parametrization is geodesic. We count images: reversal of orientation,
change of base point, and iteration do not produce new geometric
geodesics. Constant curves are excluded. Let $\Nscg(K,g)$ denote
their number, possibly infinite.

The main result of the paper is the following.

\begin{theorem}\label{thm:main}
Every smooth Riemannian metric on $K$ has at least four geometrically
distinct simple closed geodesics. Moreover, for every
$0<\varepsilon<1/4$, the metric
\begin{equation}\label{eq:metric}
 g^\varepsilon=r(x,y)(dx^2+dy^2),
 \qquad
 r(x,y)=1+\varepsilon(\cos x+\cos y),
\end{equation}
descends to $K$ and satisfies
\[
 \Nscg(K,g^\varepsilon)=4.
\]
Its four simple closed geodesics are the images of
$y=0$, $y=\pi$, $x=0$, and $x=\pi$.
\end{theorem}

\begin{corollary}\label{cor:sharp}
The optimal universal lower bound is
\[
 \min_g\Nscg(K,g)=4,
\]
where the minimum is taken over all smooth Riemannian metrics on $K$.
\end{corollary}

Section~\ref{sec:correction} identifies the issue in the published
argument and proves the lower bound. The remaining sections establish
that the metric \eqref{eq:metric} has exactly four simple closed
geodesics. The proof of this explicit count is independent of the
universal lower-bound argument.

\medskip

\noindent\textbf{Acknowledgments.}
We are very grateful to Lucas Ambrozio for suggesting the problem of determining whether the lower bound in Ballmann's work for the Klein bottle case is sharp. GF is supported by Coordenação de Aperfeiçoamento de Pessoal de Nível Superior (CAPES) during her doctoral studies. MM is supported by the S\~ao Paulo Research Foundation (FAPESP)
under grant No.~2025/21147-7. 

\noindent\textbf{AI Disclosure.}
The authors used ChatGPT-6.1 Sol (OpenAI) to assist with proofreading and
improving the English exposition of this manuscript, translating
Ballmann's article \cite{Ballmann1978} from German into English,
and exploring ideas related to the arguments presented in this
paper. All AI-assisted material was reviewed by the authors,
who take full responsibility for the final text and its
mathematical content.

\section{Ballmann's variational framework}
\label{sec:ballmann-background}

Let $M$ be a compact Riemannian surface and set $S^1=\R/\Z$.
Following Ballmann, let $PM$ denote the space of piecewise $C^1$
maps $\gamma:S^1\to M$, with the uniform topology inherited from
continuous maps. Use the energy normalization
\begin{equation}\label{eq:bg-energy}
 E(\gamma)=\int_0^1|\dot\gamma(t)|_g^2\,dt.
\end{equation}
The group $O(2)$ acts on the parameter circle by rotations and
reflection:
\[
 \gamma(t)\longmapsto\gamma(t+\theta),
 \qquad
 \gamma(t)\longmapsto\gamma(-t).
\]
The energy is invariant under these operations, so it is defined on the quotient space
\begin{equation}\label{eq:bg-quotient}
 \Sigma M=PM/O(2),
 \qquad
 q:PM\longrightarrow\Sigma M.
\end{equation}

The purpose of the quotient is to avoid counting changes of the
starting point or of the direction as new geodesics. For a prime
closed geodesic with its constant-speed parametrization, its $O(2)$
orbit determines its geometric image. For arbitrary curves the
quotient is more limited: it does not identify every possible
reparametrization. Nor does it identify a geodesic with its iterates.
In fact, for $\gamma^m(t)=\gamma(mt)$,
\[
 E(\gamma^m)=m^2E(\gamma).
\]
This distinction matters when interpreting different critical levels.

There are two different types of loops in the argument:
\[
 \underbrace{t\longmapsto\gamma(t)\in M}
 _{\text{a curve in the surface}},
 \qquad
 \underbrace{s\longmapsto[\gamma_s]\in\Sigma M}
 _{\text{a family of curves}}.
\]
A class in $H_1(M)$ concerns the first object. The $1$-cycles used
in the multiplicity argument concern the second. All curves in a
connected family can belong to the same free homotopy class in $M$,
while the family itself represents a nonzero homology class in
$\Sigma M$.

For \(\lambda>0\), set
\[
P^\lambda M=\{c\in PM:E(c)<\lambda\},
\qquad
\Sigma^\lambda M=P^\lambda M/O(2),
\]
and let \(p:P^\lambda M\to\Sigma^\lambda M\) be the quotient
projection.

\subsection{Projected homology}

Write $\mathbb F_2=\Z/2\Z$. A singular simplex
$\sigma:\Delta^k\to\Sigma M$ is called liftable\footnote{Observe that this property would always hold if \(q\) were a covering map, which need not be the case in general.} if there is a
continuous map $\widetilde\sigma:\Delta^k\to PM$ such that
\[
 q\circ\widetilde\sigma=\sigma.
\]
The lift chooses parametrized representatives continuously over the
whole simplex. Define
\begin{equation}\label{eq:bg-projected-complex}
\begin{split}
 C_k^{\mathrm{pr}}(\Sigma M;\mathbb F_2)
 =
 \operatorname{span}_{\mathbb F_2}
 \bigl\{
 \sigma:\Delta^k\to\Sigma M:
 \sigma\text{ is liftable}
 \bigr\}.
\end{split}
\end{equation}
Every face of a liftable simplex is liftable, by restricting its
lift. Hence these groups form a chain complex with the usual singular
boundary operator. Its homology is
\begin{equation}\label{eq:bg-projected-homology}
 H_k^{\mathrm{pr}}(\Sigma M;\mathbb F_2)
 =
 \frac{
 \ker(\partial:C_k^{\mathrm{pr}}\to C_{k-1}^{\mathrm{pr}})
 }{
 \operatorname{im}
 (\partial:C_{k+1}^{\mathrm{pr}}\to C_k^{\mathrm{pr}})
 }.
\end{equation}
This is the projected complex specified in
\cite[p.~42]{Ballmann1978}.

Equivalently, at the chain level,
\[
 C_*^{\mathrm{pr}}(\Sigma M;\mathbb F_2)
 =
 \operatorname{im}
 \bigl(
 q_\#:C_*(PM;\mathbb F_2)
 \longrightarrow C_*(\Sigma M;\mathbb F_2)
 \bigr).
\]
It does not follow that its homology is the image of the induced
map on homology. A chain upstairs may have nonzero boundary even
though its projected chain is a cycle. Nor should projected homology
be identified with ordinary singular homology of the quotient
without an additional argument.

The coefficients $\mathbb F_2$ remove the need to orient simplices
consistently; they do not remove the requirement that simplices
admit lifts.

\begin{example}
The horizontal family on our Klein bottle is
\[
 h_s(t)=\pi_K(2\pi t,s),
 \qquad 0\leq s\leq2\pi.
\]
Its endpoints satisfy $h_{2\pi}(t)=h_0(-t)$. Thus the displayed
path in $PM$ does not close with the same parametrization, but its
projection closes in $\Sigma M$. Its boundary vanishes after
projection, and it supplies a projected $1$-cycle. This observation establishes that it is a cycle; proving that its
homology class is nonzero is a separate step; see Proposition \ref{prop:lower}.   
\end{example}

A useful test is a continuous map
$\Phi:\Sigma M\to S^1$, or such a map on the relevant loop-space
component. Suppose a projected $1$-cycle $u$ satisfies
\[
 \Phi_*[u]\ne0
 \quad\text{in }H_1(S^1;\mathbb F_2).
\]
Then $u$ cannot be a projected boundary. Indeed, if
$u=\partial W$ with $W$ a projected $2$-chain, then
\[
 \Phi_\#u=\partial(\Phi_\#W)
\]
would be an ordinary boundary in $S^1$.

For a family parametrized by a circle, a map of odd degree therefore
detects a nonzero projected class. This is the role of the map
$\Phi$ in \eqref{eq:detector} below.

There is also a precise way to prove vanishing. Suppose a homotopy
$F_s:PM\to PM$ is $O(2)$-equivariant and descends to
$\overline F_s:\Sigma M\to\Sigma M$. For a lifted simplex
$\widetilde\sigma$, the map
\[
 (s,z)\longmapsto F_s(\widetilde\sigma(z))
\]
lifts the homotopy of its projection. Triangulating
$[0,1]\times\Delta^k$ gives the usual prism chain homotopy entirely
within the projected complex.

Thus an equivariant contraction of this kind forces positive-degree
projected homology to vanish. Contractibility of the underlying
quotient alone does not justify this conclusion unless compatibility
with the projected complex has been checked. This is why the
compatibility verification in
Proposition~\ref{lem:one-sided-space} is included.

\subsection{Curve-shortening and min-max values}

In this section, we define the min-max values for families of curves and show that its positivity implies existence of simple or twice-traversed simple closed geodesic.

\begin{definition}
Let $X$ be topological space. A continuous map $f:X\longrightarrow\Sigma M$ is called \emph{admissible} if every point \(x_0\in X\) has
an open neighborhood \(U\) on which \(f\) admits a continuous
lift
\[
\widetilde f_U:U\longrightarrow PM,
\qquad
q\circ\widetilde f_U=f|_U,
\]
satisfying the following conditions. Writing $F_U(x,s)=\widetilde f_U(x)(s),$ we require the following:
\begin{enumerate}
    \item for every \(x\in U\), the map
    \(F_U(x,\cdot):\mathbb{S}^1\to M\) is either constant
    or injective;

    \item there exists a finite subdivision $0=s_0<s_1<\cdots<s_N=1,$
    independent of \(x\in U\), such that
    \(F_U(x,\cdot)|_{[s_j,s_{j+1}]}\) is \(C^1\)
    for every \(x\in U\) and every \(j=0,\ldots,N-1\);

    \item for every \(j=0,\ldots,N-1\), the map
    \[
    U\times[s_j,s_{j+1}]\longrightarrow TM,
    \qquad
    (x,s)\longmapsto \partial_s F_U(x,s),
    \]
    is continuous.
\end{enumerate}
\end{definition}

For an admissible family \(f\), let $H_f:X\times[0,1]\longrightarrow\Sigma M$
be a homotopy obtained from Ballmann's localized curve-shortening construction (see \cite{Lusternik1966,Ballmann1978}), with \(H_f(x,0)=f(x)\). Once the choices involved in this construction have been fixed, we denote its final admissible family by
\[
D(f):X\longrightarrow\Sigma M,
\qquad
D(f)(x)=H_f(x,1).
\]
If \(f(x)\) is a constant curve, then
\(H_f(x,t)=f(x)\) for every \(t\in[0,1]\).

\begin{theorem}[See \cite{Ballmann1978}] \label{deformationtheo}
    Let $0 < l< \lambda$, and let $G(l)$ be the subset of $\Sigma M$ consisting of 
\begin{itemize}
    \item {i)} the simple closed geodesics of energy $l$, if $M$ is orientable; or
    \item {ii)} the simple closed geodesics, and the twice-traversed simple closed geodesics that can be approximated by simple closed curves, of energy $l$, if $M$ is nonorientable.
\end{itemize}

Let $U$ be neighborhood a of $G(l)$ (if $G(l) = \emptyset $, one may take $U=\emptyset$). Then there exists $\varepsilon>0$ such that, for every admissible map $f$ as above with $im f \subset \Sigma^{l+ \varepsilon}M$, one has
$$\mathrm{im} D(f) \subset U \cup \Sigma^{l- \varepsilon}M.$$
\end{theorem}

\begin{remark}
The deformation theorem also applies within a fixed
unoriented free homotopy class. Indeed, the corresponding
subset $\Sigma_a\subset\Sigma M$ is open and closed in the
uniform topology and is preserved by the deformation
homotopies. Thus the theorem and its min--max existence
corollary remain valid with $\Sigma M$ replaced by
$\Sigma_a$ and $G(l)$ replaced by $G(l)\cap\Sigma_a$.
\end{remark}

\begin{definition}
    Let $f:X\to\Sigma M$ be an admissible family with bounded energy,
and let $\nu$ be a homology class in its parameter space. Consider
cycles representing $\nu$ and admissible deformations of $f$.
After the subdivisions needed to obtain projected chains, write
$\mathcal F_f(\nu)$ for the resulting collection of image cycles.

The associated number is
\begin{equation}\label{eq:bg-minmax}
l(\nu)
 =
 \inf_{u\in\mathcal F_f(\nu)}
 \ \sup_{z\in|u|}E(z),
\end{equation}
where $|u|$ denotes the support of the cycle.
\end{definition}

For completeness, we provide a proof of the following corollary, as it was stated without proof in \cite{Ballmann1978}.

\begin{corollary}\label{cor: Existence of closed geodesic}
    If $l(\nu)>0$, then $G(l(\nu))$ is nonempty.
\end{corollary}

\begin{proof}
    Assume that $l(\nu)>0$. By the definition of $l(\nu)$,
for every \(\varepsilon>0\) there exists \(u\in \mathcal{F}_f(\nu)\)
such that
\[
l(\nu)\leq \max_{c\in |u|}E(c)<l(\nu)+\varepsilon.
\]

Suppose, for a contradiction, that \(G(l(\nu))=\varnothing\).
We may then take \(U=\varnothing\). By Theorem \ref{deformationtheo}, there
exists \(\varepsilon>0\) such that every admissible map \(f\)
satisfying $im f \subset \Sigma^{l+ \varepsilon}M$ its deformation satisfies

$$
\operatorname{im}D(f)
\subset U\cup\Sigma^{l(\nu)-\varepsilon}M
=\Sigma^{l(\nu)-\varepsilon}M.
$$

Recall that
\[
\mathcal F_f(\nu)=
\left\{
u\in C_k^{\mathrm{pr}}(\Sigma^\lambda M;\mathbb{Z}_2)
\;\middle|\;
\begin{array}{l}
u=g_\#(\widetilde{\nu}),\quad \widetilde{\nu} \ \ \text{is homologous to}\ \ \nu\\
g
\text{ is admissible and homotopic to f }
\end{array}
\right\}.
\]

Applying the deformation to the cycle $u$, we obtain the cycle $D(f)_{\#}(u)$. As established above, its support satisfies
\[
|D(f)_{\#}(u)|\subset\Sigma^{l(\nu)-\varepsilon}M.
\]

Since \(D(f)\) is admissible, we have \(D(f)_{\#}(u)\in \mathcal F_f(\nu)\).
Therefore, by the definition of \(l(\nu)\),
\[
l(\nu)\leq\max_{c\in|D(f)(u)|}E(c)<l(\nu)-\varepsilon,
\]
which is a contradiction.
\end{proof}

For an energy value $c>0$, let $\mathcal G(c)$ denote the relevant
geodesic orbits of energy $c$. Intuitively, the Theorem \ref{deformationtheo} says that for every neighborhood
$U$ of this set there is $\delta>0$ such that an admissible family
with energy below $c+\delta$ can be deformed into
\begin{equation}\label{eq:bg-deformation}
 U\ \cup\ \{E<c-\delta\}.
\end{equation}
Away from the critical set, the entire family can therefore be
pushed below the level $c$.

\subsection{Min-max inequality for subordinate classes}\label{sectionMIN-MAX INEQ}

Throughout this section, we use the cap product in the projected homology. Notice that every face of a liftable
simplex it is also liftable so the usual cap product restricts to this complex.

Let $u_1$ be a projected homology class and let
$\zeta\in H^d(\Sigma M;\mathbb F_2)$, where $d>0$. If
\begin{equation}\label{eq:bg-subordination}
 u_0=u_1\frown\zeta\ne0,
\end{equation}
we say that $u_0$ is subordinate to $u_1$. The significance of
this relation for min--max theory is that the cap product
does not enlarge the support of a chain. Consequently, a
competitor for the higher-dimensional class produces a
competitor for the subordinate class without increasing
the maximal energy.

To apply this observation to the admissible min--max
construction, we specify the compatibility of the two
classes. Let $f:X\to\Sigma M$ be an admissible family with
bounded energy, let $\nu_1\in H_k(X;\mathbb F_2)$, and set
\[
 \nu_0=\nu_1\frown f^*\zeta.
\]
Assume that $f_*\nu_0\ne0$, and write $ u_i=f_*\nu_i,$ $ c_i=l(\nu_i)$ for $i\in\{0,1\}$, where both min--max values use the same collection of
admissible deformations of $f$. Naturality of the cap product
then gives \eqref{eq:bg-subordination}.

\begin{proposition}\label{prop:bg-subordination-inequality}
With these choices, $c_0\leq c_1$.
\end{proposition}

\begin{proof}
Choose a cocycle $\xi$ representing $\zeta$, and let $v_1$
be any competitor in the construction of $c_1$. After the
subdivisions required to obtain projected chains, we may
write $v_1=F_\#a$, where $F:X\to\Sigma M$ is an admissible deformation of $f$
and $a$ is a cycle representing $\nu_1$.

The chain
\[
 a_0=a\frown F^*\xi
\]
is a cycle by the boundary formula for the cap product.
Since $F$ is homotopic to $f$, its homology class is
\[
 [a_0]
 =\nu_1\frown F^*\zeta
 =\nu_1\frown f^*\zeta
 =\nu_0.
\]
Thus $v_0:=F_\#a_0$ is a competitor in the construction
of $c_0$. At the chain level, naturality gives
\[
 v_0
 =F_\#(a\frown F^*\xi)
 =(F_\#a)\frown\xi
 =v_1\frown\xi.
\]
The cap-product formula expresses $v_0$ as a sum of faces
of simplices occurring in $v_1$. These faces are liftable,
and their images lie in $|v_1|$. Therefore $|v_0|\subseteq|v_1|$, and hence
\[
 c_0
 \leq \sup_{z\in|v_0|}E(z)
 \leq \sup_{z\in|v_1|}E(z).
\]
Taking the infimum over all competitors $v_1$ proves
$c_0\leq c_1$.
\end{proof}

We use the following lemma in the Proposition \ref{prop:lower} to find a simple closed geodesic associated to a constant admissible family.

\begin{lemma}[A lower bound for the point-class min--max]
Let $(M,g)$ be a Riemannian manifold, and let $\gamma\in PM$
represent a nontrivial free homotopy class $a$. Set
\[
u_0=[q(\gamma)]
\in H_0^{\mathrm{pr}}(\Sigma M;\mathbb F_2),
\qquad
\ell_g(a):=\inf\{L_g(\beta):\beta\text{ represents }a\}.
\]
Every projected $0$-cycle $z$ representing $u_0$ satisfies
\[
\sup_{p\in|z|}E(p)\ge\ell_g(a)^2.
\]
Consequently, any min--max value defined using projected
representatives of $u_0$ satisfies the inequality $l(\nu_0)\ge\ell_g(a)^2 \textrm{ with } f_*\nu_0=u_0.$
\end{lemma}

\begin{proof}
Define a linear map $\chi_a:C_0^{\mathrm{pr}}(\Sigma M;\mathbb F_2)
\longrightarrow\mathbb F_2$, by setting
\[
\chi_a(q(\beta))=
\begin{cases}
1,&\text{if $\beta$ represents $a$ or $a^{-1}$},\\
0,&\text{otherwise}.
\end{cases}
\]
Here $a^{-1}$ denotes the class obtained by reversing the
orientation of a representative of $a$. Thus $\chi_a$ is
well defined on the quotient by $O(2)$.

Every liftable $1$-simplex $\sigma$ admits a lift
$\widetilde\sigma:[0,1]\to PM$, which gives a free homotopy
between its endpoint curves. Consequently,
\[
\chi_a(\partial\sigma)
=\chi_a(\sigma(1))+\chi_a(\sigma(0))=0
\ \ \text{in }\mathbb F_2.
\]
Hence $\chi_a$ vanishes on projected boundaries.

If $z$ represents $u_0$, we therefore have
\[
\chi_a(z)=\chi_a(q(\gamma))=1.
\]
It follows that $|z|$ contains a point $q(\beta)$ with
$\beta$ representing $a$ or $a^{-1}$.
Since reversal preserves length,
$\ell_g(a^{-1})=\ell_g(a)$. Cauchy--Schwarz inequality gives
\[
\sup_{p\in|z|}E(p)
\ge E(\beta)
\ge\left(\int_0^1|\dot\beta(t)|_g\,dt\right)^2
=L_g(\beta)^2
\ge\ell_g(a)^2.
\]
Taking the infimum over the admissible competitors proves
the min--max inequality.
\end{proof}

\begin{remark}[Positivity on closed manifolds]
If $M$ is closed, then the injective radius $i_g:=\operatorname{inj}(M,g)>0$.
Every noncontractible loop has length at least $2i_g$.
Thus, for every nontrivial free homotopy class $a$,
\[
\ell_g(a)\ge2i_g,
\qquad
l_f(\nu_0)\ge\ell_g(a)^2\ge4i_g^2>0.
\]
\end{remark}

\section{Correction of the Klein-bottle argument}
\label{sec:correction}

We use the deformation theorem for families of embedded curves from
\cite[Theorem~1.3 and Corollary~1.4]{Ballmann1978}, including its allowance
for double covers of one-sided geodesics as limits. Throughout this
section, homology has coefficients in $\mathbb{F}_2$, and the quotient by
$O(2)$ identifies rotations and reversal of the parameter circle.

\subsection{Contractibility of the curve space in a one-side class}\label{sec: correctquotation}

In this section, we identify the flaw in Ballmann’s argument. More precisely, we show that the curve space associated with a primitive one-sided class is contractible, so the proposed homological construction cannot yield two geometrically distinct geodesics in that class.

\begin{definition}\label{def: one-two side curves}
Let $M$ be a smooth surface, let
$c \colon \mathbb{S}^1 \to \operatorname{int}(M)$
be a smooth embedding, and set $C = c(\mathbb{S}^1)$.
Let $V$ be a sufficiently small compact tubular neighborhood of $C$.
We say that:
\begin{itemize}
    \item $C$ is \emph{two-sided} if the pair $(V,C)$ is
    diffeomorphic to $  \bigl(\mathbb{S}^1 \times [-1,1], \mathbb{S}^1 \times \{0\}\bigr);$

    \item $C$ is \emph{one-sided} if the pair $(V,C)$ is
    diffeomorphic to the pair $(\mathcal{M},C_0)$ consisting
    of the Möbius band and its central curve, where
    \[
        \mathcal{M}
        = \bigl([0,1] \times [-1,1]\bigr)/{\sim},
        \qquad (0,t) \sim (1,-t),
    \]
    and $C_0 = \{[s,0] : s \in [0,1]\}$
    
\end{itemize}
\end{definition}

\textbf{Correction:} In the proof of \cite[Theorem~2.3, p.~46]{Ballmann1978}, immediately before
the final sentence, Ballmann writes:
\begin{quote}
``$b$ und $ab$ liefern wieder 1-Zykel (wie beim Torus), so da\ss{}
$b$ und $ab$ zu diesen subordiniert sind.''
\end{quote}

The asserted subordination is used to obtain two geodesics from each
of the classes denoted by $b$ and $ab$ in the quotation, and hence a
total of five. With our edge identifications, these are our classes
$a$ and $ab$, respectively. The obstruction is in the one-sided
class $ab$.

More intrinsically, among the three essential embedded classes that
are not proper powers, two are one-sided and only one is two-sided;
thus two distinct such classes cannot both support the torus
construction. The following direct computation shows that a primitive
one-sided class has no nontrivial $1$-cycle of the required kind.

\begin{proposition}\label{lem:one-sided-space}
For either primitive one-sided class $b$ or $ab$, the corresponding
space of piecewise $C^1$ loops modulo $O(2)$ is contractible.
The contraction can be performed on parametrized loops before taking
the quotient, so its positive-degree homology also vanishes for the
projected complex used in \cite{Ballmann1978}.
\end{proposition}

\begin{proof}
Let $\delta=0$ or $\pi$ and put
\[
 h_\delta(x,y)=(2\delta-x,y+2\pi).
\]
These are respectively $B$ and $AB$. Orient a loop in the relevant
unoriented class so that its vertical displacement is positive.
This selects one of its two orientations. Choose a lift
$\widetilde\gamma(t)=(x(t),y(t))$, extended each component to $\R$, satisfying
\[
 x(t+1)=2\delta-x(t),
 \qquad
 y(t+1)=y(t)+2\pi.
\]
Such a lift exists because the closing transformation is conjugate
to $h_\delta$. Set
\[
 c(\widetilde\gamma)
 =
 \int_0^1\bigl(y(t)-2\pi t\bigr)\,dt
\]
and, for $0\leq s\leq1$, define
\begin{equation}\label{eq:contraction}
 \mathcal H_s(\widetilde\gamma)(t)
 =
 \left(
 (1-s)x(t)+s\delta,\,
 (1-s)y(t)+s\bigl(2\pi t+c(\widetilde\gamma)\bigr)
 \right).
\end{equation}
The same endpoint relations hold throughout this homotopy.

Two choices of lift with closing transformation $h_\delta$ differ
by an element of its centralizer in $\Gamma$, which is
$\langle h_\delta\rangle$. Formula~\eqref{eq:contraction} commutes with this
action. It also commutes with rotations of the parameter: if
$\widetilde\gamma_\theta(t)=\widetilde\gamma(t+\theta)$, periodicity
of $y(t)-2\pi t$ gives
\[
 c(\widetilde\gamma_\theta)
 =
 c(\widetilde\gamma)+2\pi\theta.
\]
Thus the homotopy descends to the space of unoriented loops modulo
parameter rotations. At $s=1$ every loop is the image of the same
straight vertical geodesic $x=\delta$, with the same constant-speed
parametrization up to a rotation. This gives the contraction.

The homotopy on parametrized loops supplies lifts of the homotopies
of lifted simplices. Its induced chain homotopy is therefore valid
in the projected complex as well. Embeddedness need not be preserved:
the homology used in the disputed step is that of the ambient
quotient loop space.
\end{proof}

Consequently the asserted essential $1$-cycle in the class $ab$
cannot provide the additional geodesic. The flat metric gives the
same obstruction at the level of critical points; see Example \ref{example: flatmetric}.

There is also an embedded two-sided class $b^2$. The vertical curves
$x=c$, $0<c<\pi$, represent it and close after a displacement of
$4\pi$ in $y$. As $c$ approaches $0$ or $\pi$, this interval of
parametrized curves converges to double traversals of the one-sided
curves. Those are allowed limits in the deformation theorem and do
not give an additional geometric geodesic. In particular, this family
does not supply a replacement for the claimed cycle in $ab$.

\subsection{A circle of curves in the two-sided class}

Finally, we show that every Riemannian metric on the Klein bottle admits at least four geometrically distinct simple closed geodesics: at least two in the two-sided class $a$ and at least one in each of the one-sided classes $b$ and $ab$. This establishes the universal lower bound asserted in Theorem~\ref{thm:main}.

\begin{proposition}\label{prop:lower}
Every smooth Riemannian metric on $K$ has at least four geometrically
distinct simple closed geodesics.
\end{proposition}

\begin{proof}
The classes $b$ and $ab$, represented by $B$ and $AB$, contain
embedded circles. Minimization within each embedded class, using
the deformation Theorem \ref{deformationtheo}, gives a
simple closed geodesic. Neither class can degenerate to a point or
to a double cover: both are nontrivial and orientation reversing.

For the two-sided class $a=[A]$, consider the family of embedded curves
\begin{equation}\label{eq:horizontal-family}
 h_s(t)=\pi_K(2\pi t,s),
 \qquad
 t\in\R/\Z,
 \quad
 0\leq s\leq2\pi.
\end{equation}
Since $h_{2\pi}(t)=h_0(-t)$, this is a closed family after quotienting
by the $O(2)$ action on the parameter circle. The curves may be
parametrized with constant speed for the given metric.

Here is an explicit check that the resulting $1$-cycle is nontrivial.
Let $\mathcal L_a$ be the space of parametrized loops in the free
homotopy class $a$, including reversal of orientation, and set
$\Sigma_a=\mathcal L_a/O(2)$. A lift $(x(t),y(t))$ of such a loop
satisfies $y(t+1)=y(t)$. Define
\begin{equation}\label{eq:detector}
 \Phi:\Sigma_a\longrightarrow S^1,
 \qquad
 \Phi([\gamma])
 =
 \exp\left(i\int_0^1 y(t)\,dt\right).
\end{equation}
Changing the lift changes $y$ by an integral multiple of $2\pi$.
Parameter rotations and reversal leave the integral unchanged.
Consequently $\Phi$ is well defined and continuous. On
\eqref{eq:horizontal-family}, it is exactly $\Phi([h_s])=e^{is}$.
Thus this family detects the generator of $H^1(S^1;\mathbb{F}_2)$.
It also defines a cycle in the projected complex, since each
parameter interval admits the displayed parametrized lift.

Let $u_1$ be its homology class and let $\eta$ be the generator of
$H^1(S^1;\mathbb{F}_2)$. Then $u_0=u_1\frown\Phi^*\eta$ is the nonzero point
class. By Proposition \ref{prop:bg-subordination-inequality}, the min-max value for the point class $c_0$ is less than or equal to the minx-max value $c_1$ for the class $u_1$.

\textbf{Case 1:} $0<c_0<c_1$.

By Corollary~\ref{cor: Existence of closed geodesic},
applied within the horizontal free homotopy component,
there is a critical geodesic orbit at each of the levels
$c_0$ and $c_1$. Since $A$ is not conjugate to a square
in $\Gamma$, neither orbit can correspond to a
twice-traversed geodesic. Both are therefore represented
by simple closed geodesics, which are geometrically
distinct because their energies differ.

\textbf{Case 2:} $0<c=c_0=c_1$.

Set $\zeta=\Phi^*\eta\in H^1(\Sigma_a;\mathbb F_2)$, so that
$u_0=u_1\frown\zeta\ne0$. Let $\mathcal G_a(c)$ denote the
critical geodesic orbits of energy $c$ in $\Sigma_a$ relevant to
the deformation theorem. As in Case~1, none of these orbits is a
twice-traversed geodesic.

We claim that every sufficiently small neighborhood $U$ of
$\mathcal G_a(c)$ in $\Sigma_a$ contains a projected $1$-cycle $v$
such that
\[
 \langle\Phi^*\eta,[v]\rangle\ne0.
\]
Suppose instead that $\zeta$ evaluates trivially on every
projected $1$-cycle supported in $U$. Its restriction to the
projected complex of $U$ is then a coboundary. Extending the
corresponding $0$-cochain to $\Sigma_a$, we may choose a cocycle
$\xi$ representing $\zeta$ that vanishes on every liftable
$1$-simplex contained in $U$.

By Theorem~\ref{deformationtheo}, applied within $\Sigma_a$,
there is $\delta>0$ such that a competitor for $c_1$ sufficiently
close to $c$ can be admissibly deformed to a projected cycle
$v_1$ satisfying
\[
 |v_1|\subset
 U\cup\{z\in\Sigma_a:E(z)<c-\delta\}.
\]
Write $v_1=F_\#\alpha$, where $F:S^1\to\Sigma_a$ is an admissible
deformation of the horizontal family and $\alpha$ represents $[S^1]$.
Admissibility implies that $E\circ F$ is continuous. Thus, after
subdivision in the parameter circle, every simplex of $v_1$ is
contained either in $U$ or in
$\{z\in\Sigma_a:E(z)<c-\delta\}$. The simplices contained in $U$
contribute nothing to $v_1\frown\xi$, and every remaining
contribution is supported below $c-\delta$. Consequently, the
cycle
\[
 v_0:=v_1\frown\xi
\]
satisfies
\[
 |v_0|\subset\{z\in\Sigma_a:E(z)<c-\delta\}.
\]
The compatibility argument in
Proposition~\ref{prop:bg-subordination-inequality} shows that
$v_0$ is a competitor for $c_0$, since
\[
 [\alpha\frown F^*\xi]
 =[S^1]\frown F^*\zeta
 =[S^1]\frown\eta
 =[\mathrm{point}].
\]
Hence
\[
 c_0\leq\sup_{z\in|v_0|}E(z)<c-\delta<c,
\]
contrary to $c_0=c$. This proves the claim.

If $\mathcal G_a(c)$ were finite, Ballmann's local result
\cite[Lemma~1.1]{Ballmann1978} would provide a neighborhood $U$
whose projected $1$-cycles are nullhomologous in the ambient
projected complex. The class $\Phi^*\eta$ would evaluate
trivially on all such cycles, contradicting the claim.
Therefore $\mathcal G_a(c)$ is infinite. Since all its orbits
are simple closed geodesics in the primitive class $a$, and a
simple geodesic image at the fixed energy $c$ determines a
single orbit modulo rotations and reversal, there are
infinitely many geometrically distinct simple closed geodesics
in this class. In particular, there are at least two.

The two geodesics in the class $a$ are two-sided, and hence differ
from the two already obtained in the one-sided classes. This gives
at least $2+1+1=4$ geometrically distinct simple closed geodesics.
\end{proof}

\begin{remark}
The class $a$ supplies \emph{at least} two geodesics by this argument;
the other two classes each supply at least one. These are lower
bounds for an arbitrary metric, not assertions of uniqueness within
the classes. The explicit metric below realizes exactly these
numbers and has no simple closed geodesic in the class $b^2$.
\end{remark}

\section{The four geodesics for $g^{\varepsilon}$}

In this section, we show that the metric $g^{\varepsilon}$, with $0<\varepsilon<1/4$, has exactly four geometrically distinct simple closed geodesics. We first identify these four geodesics and then rule out any others.

\subsection{A criterion for coordinate lines to be geodesics}
For the rest of the paper, assume $0<\varepsilon<1/4$ and use the Riemannian metric in
\eqref{eq:metric}. Set $\rho_0=1-2\varepsilon,
 \rho_1=1+2\varepsilon.$ Then $0<\rho_0\leq r\leq\rho_1$.

\begin{lemma}\label{lem:coordinate}
    For a conformal metric $g=r(x,y)(d x^2+d y^2)$, the horizontal
coordinate line $y=y_0$ is an unparametrized geodesic if and only if
\[
 \partial_y r(x,y_0)=0\quad\text{for every }x.
\]
Similarly, the vertical coordinate line $x=x_0$ is an unparametrized geodesic if
and only if
\[
 \partial_x r(x_0,y)=0 \ \ \text{for every }y.
\]
\end{lemma}

\begin{proof}

We shall focus on the horizontal coordinate line $y= y_0$, as the vertical case is entirely analogous. 
In coordinates, the metric $g = r(x,y)(dx^2 + dy^2)$ takes the form
\[
g_{ij} = e^{2u(x,y)}\delta_{ij},
\qquad
g^{ij} = e^{-2u(x,y)}\delta^{ij},
\]
where
\[
u(x,y) = \frac{1}{2}\log r(x,y),
\]
and $(g^{ij})$ denotes the inverse matrix of $(g_{ij})$.

The Christoffel symbols are given by
\[
\Gamma_{ij}^{k}
=
\frac{1}{2}
\sum_{\ell}
g^{k\ell}
\left(
\frac{\partial g_{j\ell}}{\partial x^i}
+
\frac{\partial g_{i\ell}}{\partial x^j}
-
\frac{\partial g_{ij}}{\partial x^\ell}
\right).
\]
In particular, we have $\Gamma_{xx}^{x}=\partial_x u,\ \ 
\Gamma_{xx}^{y}=-\partial_y u.$

Since we are considering unparametrized geodesics, recall that a regular curve is an unparametrized
geodesic if and only if
\[
\nabla_{\dot\gamma}\dot\gamma
=\lambda(t)\dot\gamma
\]
for some smooth function $\lambda$.

We may use the coordinate parametrization $\gamma(t)=(t,y_0)$.

Using an affine connection $\nabla$ in local coordinates, in terms of Christoffel symbols, we have
\[
\nabla_{\dot\gamma}\dot\gamma
= \nabla_{\partial_x} \partial_x = \Gamma_{xx}^x \partial_x + \Gamma_{xx}^y \partial_y =
(\partial_xu)\partial_x
-(\partial_yu)\partial_y.
\]

Since $\dot\gamma=\partial_x$, this acceleration is tangent to $\gamma$ if and only if $\partial_yu(x,y_0)=0$ for every $x$ and this implies that $\partial_y r(x,y_0)=0\text{ for every }x.$

\end{proof}

Since $r_x=-\varepsilon\sin x$ and $r_y=-\varepsilon\sin y$, the
coordinate geodesics project to exactly the following four images:
\begin{equation}\label{eq:four}
\begin{aligned}
 H_0(t)&=\pi_K(t,0),
 &H_\pi(t)&=\pi_K(t,\pi),\\
 V_0(t)&=\pi_K(0,t),
 &V_\pi(t)&=\pi_K(\pi,t),
\end{aligned}
\qquad 0\leq t\leq2\pi.
\end{equation}
For $H_0$ and $H_{\pi}$ the deck transformation correspond to $A$, for $V_0$ is $B$ and $V_{\pi}$ is $AB$. Therefore, each curve is injective before its endpoints are identified. It also shows that the four images are distinct. The horizontal curves are two-sided; the vertical curves
are one-sided.

\subsection{Separated equations and exclusion of additional geodesics}

To separate the variables and apply comparison theorems, we write the conformal factor of the metric as
\begin{equation}\label{eq:FG}
 F(x)=\tfrac12+\varepsilon\cos x,
 \qquad
 G(y)=\tfrac12+\varepsilon\cos y,
 \qquad
 r(x,y)=F(x)+G(y).
\end{equation}
For the Riemannian Klein bottle $(K,g_{\varepsilon})$, the Legendre transformation is given by
\begin{equation}
\Phi: TK \rightarrow T^*K,
\qquad
\Phi(v) = g^{\varepsilon}(v,\cdot).
\end{equation}
The kinetic-energy Lagrangian $L: TK \rightarrow \mathbb{R}$, given by
$$
L(v) = \frac{1}{2}g^{\varepsilon}(v,v),
$$
has a well-defined Legendre transform, given by the Hamiltonian
$H: T^*K \rightarrow \mathbb{R}$, defined by
\[
H(p)=\frac{1}{2}g^{\varepsilon}_*(p,p),
\]
where $g^{\varepsilon}_*$ denotes the dual metric. Explicitly,
$g^{\varepsilon}_*(p,\xi)=g^{\varepsilon}(p^{\#},\xi^{\#})$, where $p^{\#}$ is the unique tangent vector to $K$ satisfying
$p=g^{\varepsilon}(p^{\#},\cdot)$.

In coordinates $(x,y)$ on $K$, the Legendre transform takes the form
\[
\Phi(x,y,v_1,v_2) = (x,y,p_x,p_y),
\]
where $p_x=\frac{\partial L}{\partial v_1}=rv_1$ and
$p_y=\frac{\partial L}{\partial v_2}=rv_2$.

Unit-speed closed geodesics correspond bijectively to periodic orbits of the Hamiltonian flow on \(H^{-1}(1/2)\), with the correspondence preserving periods. Consider a unit-speed geodesic parametrized by arclength $t$. The Hamiltonian over the geodesic is given by
$$
H=\frac{p_x^2+p_y^2}{2r}.
$$
Hamilton's equations can be computed as
\[
\begin{cases}
\displaystyle \dot{x}=\frac{\partial H}{\partial p_x},
&
\displaystyle \dot{p}_x=-\frac{\partial H}{\partial x},
\\[6pt]
\displaystyle \dot{y}=\frac{\partial H}{\partial p_y},
&
\displaystyle \dot{p}_y=-\frac{\partial H}{\partial y},
\end{cases}
\quad\Longrightarrow\quad
\begin{cases}
\displaystyle \dot{x}=\frac{p_x}{r},
&
\displaystyle \dot{p}_x=\frac{p_x^2+p_y^2}{2r^2}F'(x),
\\[6pt]
\displaystyle \dot{y}=\frac{p_y}{r},
&
\displaystyle \dot{p}_y=\frac{p_x^2+p_y^2}{2r^2}G'(y).
\end{cases}
\]
On the energy level $H= \frac{1}{2}$, we have $p_x^2+p_y^2=r$, and the system reduces to
\[
\begin{cases}
\displaystyle \dot{x}=\frac{p_x}{r},
&
\displaystyle \dot{p}_x=\frac{F'(x)}{2r},
\\[6pt]
\displaystyle \dot{y}=\frac{p_y}{r},
&
\displaystyle \dot{p}_y=\frac{G'(y)}{2r}.
\end{cases}
\]
Moreover,
\[
p_x^2+p_y^2=F(x)+G(y)
\quad\Longrightarrow\quad
j:= p_x^2-F(x)=G(y)-p_y^2.
\]
The quantity $j$ is constant over this closed geodesic on a fixed energy level.

We now introduce a new parameter $\tau$ defined by $d\tau/dt=1/r$, thereby eliminating the factor $1/r$ from Hamilton's equations. Now, primes denote differentiation with respect to $\tau$. We obtain
\begin{align}
 (x')^2&=F(x)+j,
 &(y')^2&=G(y)-j,
 \label{eq:first-order}\\
 x''&=-\frac{\varepsilon}{2}\sin x,
 &y''&=-\frac{\varepsilon}{2}\sin y.
 \label{eq:pendula}
\end{align}
The second-order equations follow directly from Hamilton's equations and remain valid at turning points. All subsequent periods are measured with respect to $\tau$.

\begin{proposition}\label{prop:contractible}
The metric $g^{\varepsilon}$ has no nonconstant contractible closed
geodesic.
\end{proposition}

\begin{proof}
Let $\gamma$ be contractible closed geodesic and let $\tilde{\gamma}$ be a closed lift to $\mathbb{R}^2$. Since $r(x,y)>0$, the coordinate time change is positive, so both $x$ and $y$ are periodic as real-valued functions of $\tau$. Choose a maximum point where $x'=0$ and denote its $x$-coordinate by $x_*$. Separately choose a point where $y'=0$ and denote its $y$-coordinate by $y_*$.

Equation \eqref{eq:first-order} gives
$F(x_*)+j=0$ and $G(y_*)-j=0$. Hence
\[
 r(x_*,y_*)=F(x_*)+G(y_*)=0,
\]
contrary to $r\geq\rho_0>0$. The two turning points need not occur
at the same time: positivity holds at every point of $\R^2$.
\end{proof}

We now state a lemma used to estimate the periods of the separeted components of a geodesic.

\begin{lemma}\label{lem:period}
Let $z$ be a nonconstant real-valued $T$-periodic solution of
\[
 z''=\frac{\varepsilon}{2}V'(z),
\]
where $V$ is smooth and $-V''\leq M$ on the range of $z$, with
$M>0$. Then
\begin{equation}\label{eq:period-bound}
 T\geq2\pi\sqrt{\frac{2}{\varepsilon M}}.
\end{equation}
\end{lemma}

\begin{proof}
The periodic function $v=z'$ is nonzero and has mean zero. It
satisfies $v''+q(\tau)v=0$, where
\[
 q=-\frac{\varepsilon}{2}V''(z)
 \leq\frac{\varepsilon M}{2}.
\]
Integration by parts and the periodic Wirtinger inequality yield
\[
 \left(\frac{2\pi}{T}\right)^2\int_0^T v^2\,d\tau
 \leq\int_0^T(v')^2\,d\tau
 =\int_0^T qv^2\,d\tau
 \leq\frac{\varepsilon M}{2}\int_0^T v^2\,d\tau.
\]
Canceling the positive integral proves \eqref{eq:period-bound}.
\end{proof}

For either equation in \eqref{eq:pendula}, take $V(z)=\cos z$ and
$M=1$. Every nonconstant periodic real-valued coordinate therefore
has period at least $2\pi\sqrt{2/\varepsilon}$.

\subsection{Exclusion of additional essential simple closed geodesics}

The orientation double cover is
\begin{equation}\label{eq:cover}
 p:\Khat=\R^2/\langle A,B^2\rangle
   =\R^2/(2\pi\Z\times4\pi\Z)
   \longrightarrow K,
\end{equation}
with deck involution $\sigma([x,y])=[-x,y+2\pi]$.
We use $(1,0)$ and $(0,1)$ for its horizontal and vertical homology
generators. In particular, the vertical generator has displacement
$4\pi$.

\begin{lemma}\label{lem:topology}
If $\gamma$ is an essential embedded circle on $K$, every component
of $p^{-1}(\gamma)$ is an essential embedded circle on $\Khat$ with
primitive class $(\pm1,0)$ or $(0,\pm1)$. If $\gamma$ is one-sided,
its preimage is a single circle, double-covering $\gamma$, and its
class is $(0,\pm1)$.
\end{lemma}

\begin{proof}
For a two-sided circle the preimage has two components, each mapping
homeomorphically to $\gamma$; for a one-sided circle it is connected
and has covering degree two. The components are embedded and
essential. Essentiality follows from injectivity on fundamental
groups and, in the one-sided case, from the fact that the square of
an orientation-reversing deck transformation is nontrivial.

An essential embedded circle on a torus has primitive homology class
$(m,n)\ne(0,0)$. The involution induces
$\sigma_*(m,n)=(-m,n)$. If there are two components, $C$ and
$\sigma(C)$, they are disjoint, and therefore
\[
 0=[C]\mathbin{\cdot}[\sigma(C)]
   =
   \det
   \begin{pmatrix}
   m&-m\\
   n&n
   \end{pmatrix}
   =2mn.
\]
Primitivity gives the stated alternatives. If the preimage is
connected, the restriction of $\sigma$ is a free involution of a
circle and hence preserves its orientation. Thus
$(-m,n)=(m,n)$, and primitivity gives $(m,n)=(0,\pm1)$.
\end{proof}

In particular, the possible free homotopy classes in $K$ with essential embedded circles are only $a$, $b$, $ab$ and $b^2$. The distinction is particularly relevant because each class imposes different periodicity conditions on the coordinates of a lifted geodesic. We will use these conditions, together with the period estimates established in Lemma \ref{lem:period}, to classify all simple closed geodesics of the metric \(g^\varepsilon\).

The following proposition completes the proof of Theorem~\ref{thm:main}, establishing that the optimal universal lower bound for the number of geometrically distinct simple closed geodesics on a Riemannian Klein bottle is four.

\begin{proposition}\label{prop:classification}
Let $K$ be the Klein bottle equipped with the metric
\[
g^\varepsilon
=
\bigl(1+\varepsilon(\cos x+\cos y)\bigr)
(dx^2+dy^2),
\qquad 0<\varepsilon<\frac14.
\]
Then the only geometrically distinct simple closed
geodesics of $g^\varepsilon$ are $H_0, H_\pi, V_0 \textrm{ and } V_\pi$.
\end{proposition}

\begin{proof}

By Proposition \ref{prop:contractible}, there are no nonconstant
contractible closed geodesics. It therefore
suffices to consider the essential simple
free homotopy classes $a,b,ab$ and $b^2$.

\medskip
\noindent

\textit{Case 1: The primitive horizontal class $a$}

Let $\gamma$ be a simple closed geodesic in this class. Choose an orientation and a lift such that, for some $T>0$,
\[
x(\tau+T)=x(\tau)+2\pi,
\qquad
y(\tau+T)=y(\tau).
\]

Suppose that $y$ is nonconstant. Since $y$ is periodic, there exists a turning value $y_*$ at which $y'=0$. The equation \ref{eq:first-order} give $j=G(y_*)$, and hence $(x')^2 = F(x)+G(y_*)=r(x,y_*) \geq\rho_0$.

Thus $x'$ never vanishes. Since $x$ increases by $2\pi$ over one period, we have $x'>0$.
Consequently,
\[
2\pi
=
\int_0^T x'\,d\tau
\geq T\sqrt{\rho_0},
\]
which yields
\[
T\leq\frac{2\pi}{\sqrt{\rho_0}}.
\]

On the other hand, applying Lemma \ref{lem:period} to the nonconstant periodic function $y$ with $V(z)=\cos(z)$ and $M=1$, gives
\[
T\geq2\pi\sqrt{\frac2\varepsilon}.
\]

These inequalities are incompatible with $\varepsilon<2\rho_0$. Therefore $y$ must be constant. The equation for $y$ then the equation \ref{eq:pendula} implies $\sin y=0$, so that
\[
y=0\quad\text{or}\quad y=\pi
\pmod{2\pi}.
\]
Thus the only geodesic images in the primitive horizontal class are $H_0$ and $H_\pi$.

\medskip
\noindent

\textit{Case 2: The primitive one-sided classes $b$ and $ab$}

Let $\gamma$ be a simple closed geodesic in one of the primitive one-sided classes. Choose a lift and an orientation such that
\[
x(\tau+T)=-x(\tau)+2\pi m,
\qquad
y(\tau+T)=y(\tau)+2\pi,
\]
where $m\in\{0,1\}$.

Applying the first identity twice, we obtain $x(\tau+2T)=x(\tau)$. Thus $x$ is $2T$-periodic.

Suppose that $x$ is nonconstant. Choose a turning value $x_*$, so that $x'=0$ at some point with $x=x_*$. The equation \ref{eq:first-order} implies $j=-F(x_*).$ As in the previous case, $(y')^2 \geq\rho_0$. Since $y$ increases by $2\pi$, we have $y'>0$. Hence
\[
T\leq\frac{2\pi}{\sqrt{\rho_0}}.
\]
However, Lemma \ref{lem:period} applied to the non constant $2T$-periodic function $x$ gives
\[
T\geq\pi\sqrt{\frac2\varepsilon}.
\]
Since $2\varepsilon<\rho_0$, we have
\[
\frac{2\pi}{\sqrt{\rho_0}}
<
\pi\sqrt{\frac2\varepsilon},
\]
a contradiction. Therefore $x$ is constant. The equation \ref{eq:pendula} implies $\sin x=0$, so

\[
x=0\quad\text{or}\quad x=\pi
\pmod{2\pi}.
\]

These two possibilities correspond precisely to the geodesics $V_0$ and $V_\pi$, representing the two primitive one-sided classes.

\textit{Case 3: The two-sided vertical class $b^2$}

As in the cases of the one-sided classes $b$ and $ab$, the preceding argument shows that any closed geodesic in the two-sided vertical class must have the same geometric image as either $V_0$ or $V_\pi$. However, both curves close after a vertical displacement of $2\pi$, whereas a geodesic representing the two-sided vertical class must have a vertical displacement of $4\pi$. Consequently, such a geodesic would traverse $V_0$ or $V_\pi$ twice and therefore would not be simple. Hence, the two-sided vertical class contains no simple closed geodesics.

\end{proof}

Proposition~\ref{prop:contractible} excludes the contractible class,
and the two arguments above exhaust the essential simple classes.
Together with the four embedded coordinate geodesics, this proves
the assertion $\Nscg(K,g^\varepsilon)=4$ in
Theorem~\ref{thm:main}. Proposition~\ref{prop:lower} proves its
universal lower-bound assertion.

Corollary~\ref{cor:sharp} follows immediately: the lower bound holds
for every metric, and every member of \eqref{eq:metric} attains it.

\bibliographystyle{plain}
\bibliography{bibliography.bib}

@article{LS1929,
  author  = {Lusternik, L. and Schnirelmann, L.},
  title   = {Sur le probl{\`e}me de trois
             g{\'e}od{\'e}siques ferm{\'e}es
             sur les surfaces de genre 0},
  journal = {Comptes Rendus de l'Acad{\'e}mie
             des Sciences de Paris},
  volume  = {189},
  pages   = {269--271},
  year    = {1929}
}

@article{hadamard,
    author = {Hadamard, J.},
    title = {Les surfaces à courbures                opposées et leurs                   lignes géodésique},
    journal = {Journal Mathematic Pure},
    pages   = {27--75},
    year = {1898}
}

@article{birkhoff,
  author  = {Birkhoff, George D.},
  title   = {Dynamical Systems with Two Degrees of Freedom},
  journal = {Transactions of the American Mathematical Society},
  volume  = {18},
  number  = {2},
  pages   = {199--300},
  year    = {1917}
}

@article{Franks1992,
  author  = {Franks, John},
  title   = {Geodesics on {$S^2$} and Periodic Points
             of Annulus Homeomorphisms},
  journal = {Inventiones Mathematicae},
  volume  = {108},
  pages   = {403--418},
  year    = {1992},
  doi     = {10.1007/BF02100612}
}

@article{Bangert1993,
  author  = {Bangert, Victor},
  title   = {On the Existence of Closed Geodesics
             on Two-Spheres},
  journal = {International Journal of Mathematics},
  volume  = {4},
  number  = {1},
  pages   = {1--10},
  year    = {1993},
  doi     = {10.1142/S0129167X93000029}
}

@article{BangertHingston1984,
  author  = {Bangert, Victor and Hingston, Nancy},
  title   = {Closed Geodesics on Manifolds
             with Infinite Abelian Fundamental Group},
  journal = {Journal of Differential Geometry},
  volume  = {19},
  number  = {2},
  pages   = {277--282},
  year    = {1984}
}

@article{Taimanov1992,
  author  = {Taimanov, I. A.},
  title   = {Closed Extremals on Two-Dimensional Manifolds},
  journal = {Russian Mathematical Surveys},
  volume  = {47},
  number  = {2},
  pages   = {163--211},
  year    = {1992}
}

@article{Ballmann1978,
  author  = {Ballmann, Werner},
  title   = {Doppelpunktfreie geschlossene Geod{\"a}tische
             auf kompakten Fl{\"a}chen},
  journal = {Mathematische Zeitschrift},
  volume  = {161},
  pages   = {41--46},
  year    = {1978},
  doi     = {10.1007/BF01175611}
}

@book{Lusternik1966,
  author    = {L. A. Lusternik},
  title     = {The Topology of the Calculus of Variations in the Large},
  series    = {Translations of Mathematical Monographs},
  volume    = {16},
  publisher = {American Mathematical Society},
  address   = {Providence, RI},
  year      = {1966}
}

\bigskip

\noindent
\textsc{Gabriella Freitas}\\
IMPA, Instituto de Matemática Pura e Aplicada.\\
\textit{E-mail:}
\href{mailto:autor2@universidade.br}
     {\texttt{gabriella.silva@impa.br}}

\medskip
     
\noindent
\textsc{Marcelo Miranda}\\
IME, Universidade de S\~ao Paulo, Brazil.\\
\textit{E-mail:}
\href{mailto:autor1@universidade.br}
     {\texttt{marcelo.miranda@ime.usp.br}}

\medskip

\end{document}